\documentclass[12pt,a4paper,preprintnumbers,showkeys]{revtex4-1} %% Other options: showpacs,twocolumn,superscriptaddress,

\usepackage{tikz}
\usepackage[fleqn]{amsmath} %% {align}
\usepackage{amssymb} %% 特殊符号, 例如 花体
\usepackage{mathrsfs} %% [MengQR 2017.4.12] \mathscr{D}
\usepackage{bm}% bold math %% 矢量符号 如 $\bm{x}$
\usepackage{graphicx} %% \includegraphics
\usepackage{subfigure}  %% 子图 Fig. 1a, 1b
\usepackage{pict2e} %% 非常规的画图命令 [DQLu 2011.1.2]
\usepackage{ulem} %% [MengQR 2017.4.12] \sout{}
\usepackage{cases}
\usepackage[pdfstartview=FitH,
            bookmarksnumbered=true,
            bookmarksopen=true,
            colorlinks=true, %注释掉此项则交叉引用为彩色边框(将colorlinks和pdfborder 同时注释掉)
            pdfborder=001,   %注释掉此项则交叉引用为彩色边框
            linkcolor=blue,  %引用链接为蓝色
            citecolor=blue,   %引用链接为蓝色
            urlcolor=blue    %网址和电邮为蓝色
            ]{hyperref}      %超链接
\usepackage[left=2.5cm, right=2.5cm, top=2.5cm, bottom=2.2cm]{geometry}
\newtheorem{definition}{Definition}[section]
\newtheorem{lemma}{Lemma}[section]
\newtheorem{theorem}{Theorem}[section]

\begin{document}

\title{On the spectral radius of uniform weighted hypergraph}

\author{Rui Sun, Wen-Huan Wang}\thanks{Corresponding author. Email: whwang@shu.edu.cn}

\affiliation{Department of Mathematics, Shanghai University, Shanghai 200444, China}

\date{\today}
\begin{abstract}
 Let $\mathbb{Q}_{k,n}$ be the set of the connected  $k$-uniform weighted hypergraphs with $n$ vertices, where $k,n\geq 3$.
  For a hypergraph $G\in \mathbb{Q}_{k,n}$, let $\mathcal{A}(G)$, $\mathcal{L} (G)$ and $\mathcal{Q} (G)$ be its adjacency tensor, Laplacian tensor  and signless Laplacian tensor, respectively.  The spectral radii of $\mathcal{A}(G)$ and $\mathcal{Q} (G)$ are investigated.
  Some basic properties of  the $H$-eigenvalue, the $H^{+}$-eigenvalue and the $H^{++}$-eigenvalue  of  $\mathcal{A}(G)$, $\mathcal{L} (G)$ and $\mathcal{Q} (G)$ are presented.  Several lower and upper bounds of the  $H$-eigenvalue, the $H^{+}$-eigenvalue and the $H^{++}$-eigenvalue for $\mathcal{A}(G)$, $\mathcal{L} (G)$ and $\mathcal{Q} (G)$   are established.  The largest $H^{+}$-eigenvalue of $\mathcal{L} (G)$ and the smallest  $H^{+}$-eigenvalue of $\mathcal{Q} (G)$ are characterized.  A relationship among the $H$-eigenvalues of $\mathcal{L} (G)$, $\mathcal{Q} (G)$ and $\mathcal{A} (G)$ is also given.
 \end{abstract}


\maketitle

%%%%%%%%%%%%%%%%%%%%%%%%%%%%%%%%%% Main Body Text %%%%%%%%%%%%%%%%%%%%%%%%%%%%%%%%%%%%%%%%%%%%%%%%%%%%%%%%%%%%%%%%%%%%%%%%%%%%

\section{Introduction}

    Weighted hypergraphs are  a natural  extension  of   hypergraphs.
    They are   of interest in real life and have many applications in graph theory.
   For example,  the circuit is mathematically modeled by a weighted hypergraph      and weighted hypergraphs are closely related to the specific application of circuit division \cite{2012-Ming-p74}.

   A weighted hypergraph is obtained from  a  hypergraph $G^{\vartriangle}=\left(V(G^{\vartriangle}),E(G^{\vartriangle})\right)$ by assigning a weight (namely, nonzero real number) to each edge of $G^{\vartriangle}$.
   We denote such a weighted hypergraph by $G=(V(G), E(G), W(G))$, where
   $V(G)=V(G^{\vartriangle})=\{v_1, \cdots,v_n\}$,
   $E(G)=E(G^{\vartriangle})=\{e_1,e_2,\cdots,e_m\}$,
   and $W(G)=\{w_G(e)\in\mathbb{R}:~e\in E(G)\}$ are  the vertex set, the edge set, and the weight set of $G$, respectively. %the set of weights  of $G$.
   Here $w_G(e)$ is the weight on the edge $e$ of  $G$ and  $\mathbb{R}$ is the set of real numbers.
    A  weighted hypergraph is  simple if it has no loops or multiple edges.
    In this paper,  we consider the weighted hypergraph which is simple and connected and satisfies that the weight of each edge is  a positive real number.

    A simple weighted hypergraph $G$ is  $k$-uniform if each edge of $G$ has $k$ vertices, where $k\geq 2 $. If $k=2 $, then $G$ is a simple weighted graph.
    \textcolor{black}{A hypergraph $G$ is called}
    linear if any two edges in $G$  intersect on at most one common vertex.
    \textcolor{black}{A hypergraph $G$ is connected}
    if there exists a path between  every pair of vertices in $V(G)$.
     Here  a path of length $p$ ($p\geq 1$) between $v_{1}$ and $v_{p+1}$ is denoted by $P=(v_{1},e_{1},v_{2},\ldots,v_{p},e_{p},v_{p+1})$, where all $v_{i}$ and all $e_{i}$ are distinct, and $v_i,v_{i+1}\in e_i$ for $1\leq i\leq p$.

   Let $G$ be a weighted hypergraph and $u,v\in V(G)$.
   \textcolor{black}{A vertex $v$ is said}
   %$v$ is said
   to be incident with an edge $e\in E(G)$ if $v\in e$.
    If  $\{u,v\}\subseteq e\in E(G)$,  then we say that  $u$ and $v$ are adjacent.
  Let $E_G(v)$ be the set of all the edges incident with $v$ of $G$, i.e., $E_G(v)=\{e\in E(G):~v\in e\}$.
  The degree  of $v$ is denoted by $d_{G}(v)$.
   Namely  $d_{G}(v)=|E_G(v)|$.
  If each vertex of $G$ has degree $r$ $(r\geq 1)$, then we say that $G$ is $r$-regular.
  We use $N_G(u)$ to denote the set of vertices which are adjacent with $u$, where $u\in V(G)$.
  For simplicity, let $\bigtriangleup = \max\limits_{v\in V(G)} |E_G (v)|$
 and $W_0 =\max\limits_{e\in E(G)} w_{G}(e).$
 Hereinafter, if each edge of $G$ has the same weight, then we denote the weight by $W_0$.
 For $v_i\in V(G)$,  let
 $w_{v_i}=\sum\limits_{e\in E_{ G} (v_i) }  w_{ G } (e)$
  and we call $w_{v_i}$  the weight of vertex $v_i$ of $G$, where $i = 1 , \ldots , n$.
 Let $\alpha =\max\{ w_{v_i}: i\in[n]\}$
and $\delta=\min\{w_{v_i}: i\in[n] \},$ where $[n]=\{1,2,\cdots, n\}$.

 A real tensor (or hypermatrix) $\mathcal{A}=(a_{i_1i_2\cdots i_k})$ of order $k$ and dimension $n$  is a multi-dimensional array
  with entries $a_{i_1i_2\cdots i_k}$, where  $a_{i_1i_2\cdots i_k}\in \mathbb{R}$ with $i_1,i_2,\cdots, i_k\in [n]$.
    The concept of tensor eigenvalues and the spectra of tensors were  introduced by   Qi \cite{2005-Qi-p1302} and Lim \cite{2005-Lim-p129} in 2005 independently.
  Let $\mathbb{C}$ be the set of  complex numbers and  $ \bm{x}=(x_1,x_2,\ldots,x_n)^{\textrm{T}}\in \mathbb{C}^n$
   an  $n$-dimensional complex column vector.
    Let  $ \bm{x}^{[k]}=(x^{k}_1,x^{k}_2,\cdots,x^{k}_n)^{\textrm{T}}$,
  where $k$ is a positive integer. Then
  $\mathcal{A}\bm{x}$
  is a vector in $\mathbb{C}^n$ whose $i$-th component is given by
  \begin{align}\label{444444}
  (\mathcal{A}\bm{x})_i=\sum^n_{i_2,\ldots,i_k=1}a_{ii_2\cdots i_k}x_{i_2}\cdots x_{i_k},~\mbox{for~each}~i\in [n].
  \end{align}
 Furthermore, we have
     \begin{equation}\label{2.2}
  \bm{x}^{\textrm{T}} (\mathcal{A}\bm{x})=\sum_{i_1,i_{2},\ldots,i_{k}=1}^{n}a_{i_1 i_{2}\ldots i_{k}}x_{i_{1}}\cdots x_{i_{k}}.
  \end{equation}

  Let $\mathbb{T}_{k,n}$ be the set of  tensors of order $k$ and dimension $n$, where $k,n\geq3$.

\begin{definition}\label{definition2.3}
 Let $\mathcal{A} \in \mathbb{T}_{k,n}$,
 {where $k,n\geq3$}. If for any vector $\bm{x} \in \mathbb{R}^{n}$,  we have
$$\bm{x}^{\textrm{T}} (\mathcal{A}\bm{x})=\sum\limits_{i_{1},i_{2},\ldots,i_{k}=1}^{n}a_{i_{1}i_{2}\ldots i_{k}}x_{i_{1}}x_{i_{2}}\cdots x_{i_{k}}
\geq 0,$$
 then $\mathcal{A}$ is called a positive semi-definite tensor.
 If for any vector $\bm{x} \in \mathbb{R}^{n}$ and $\bm{x} \neq \bm{0}$,  we have
 $\bm{x}^{\textrm{T}} (\mathcal{A}\bm{x}) > 0,$
then $\mathcal{A}$ is said to be a positive definite tensor.
\end{definition}

 For $\lambda\in \mathbb{C}$ and  $\bm{x}\in \mathbb{C}^n$, if   they satisfy  $\mathcal{A}\bm{x}=\lambda \bm{x}^{[k-1]}$, namely,  $(\mathcal{A}\bm{x})_{i}=\lambda x_i^{k-1}$ for  any $i\in[n]$,  then   $\lambda$ is called an eigenvalue of  $\mathcal{A}$ and $\bm{x}$ an
 eigenvector of  $\mathcal{A}$ corresponding to  $\lambda$.
  The largest modulus of the eigenvalues of $\mathcal{A}$ is
 \textcolor{black}{called}
 the spectral radius of $\mathcal{A}$.
 If  $\bm{x}$ is a real eigenvector of $\mathcal{A}$, then  $\lambda$ is also real and is referred to as  an $H$-eigenvalue and  $\bm{x}$ an  $H$-eigenvector.
 Let {$\mathbb{R}^{n}_{+}=\{\bm{x}\in \mathbb{R}^{n}: x_i\geq 0, i\in [n]\}$ }
 and {$\mathbb{R}^{n}_{++}=\{\bm{x}\in \mathbb{R}^{n}: x_i> 0 , i\in [n]\}$.}
  If $\bm{x}\in  \mathbb{R}^{n}_{+}$, then $\lambda$ is  an $H^{+}$-eigenvalue of $\mathcal{A}$.
  If  $\bm{x}\in  \mathbb{R}^{n}_{++}$, then $\lambda$ is  an $H^{++}$-eigenvalue of $\mathcal{A}$.

   The adjacency tensor of a $k$-uniform weighted hypergraph $G$ is defined as follows.

   \begin{definition}\label{definition2.1}%\cite{Cooper2012h}
   Let $G$ be a $k$-uniform weighted hypergraph with $n$ vertices.
   The adjacency tensor of $G$ is the $k$-ordered  and $n$-dimensional adjacency tensor $\mathcal{A}(G)=(a_{i_1i_2\cdots i_k})$ whose $(i_1i_2\cdots i_k)$-entry is
  \begin{equation}\label{A}
   a_{i_1i_2\cdots i_k}= \left\{
  \begin{array}{ll}
  \dfrac{w_{G}(e)}{(k-1)!},~~~~&\mbox{if}~ e=\{i_1,i_2,\cdots, i_k\}\in E(G), \\
  0,                   & \mbox{otherwise},
  \end{array}
  \right.
  \end{equation}
  where each $i_j$ runs from 1 to $n$ for $j\in [k]$.
   \end{definition}

   In  Definition \ref{definition2.1}, if $w_{G}(e)=1$ for each edge of $G$, then  $\mathcal{A}(G)$ is just the tensor defined by  Cooper and Dutle \cite{Cooper2012h}  in 2012 for a  $k$-uniform hypergraph with $n$ vertices.
    For a  real tensor  $\mathcal{A}=(a_{i_1i_2\cdots i_k})$, if $a_{i_1i_2\cdots i_k}$ is invariant under any permutation of the indices   $i_{1}, i_{2}, \ldots, i_{k}$, then
     $\mathcal{A}$ is said to be symmetric. A tensor is called nonnegative if all its entries are nonnegative.
     Let $G$ be a $k$-uniform weighted hypergraph.
      Obviously, the adjacency tensor $\mathcal{A}(G)$ of $G$ is always nonnegative  and symmetric.
      The spectral radius of $\mathcal{A}(G)$, denoted by  $\rho (G)$, is called the spectral radius of $G$.

  Inspired by the definitions of the  Laplacian tensor   and the signless Laplacian tensor
  \textcolor{black}{of a $k$-uniform hypergraph}
  which were introduced by Qi \cite{2014-Qi-p1045}, in this paper,
  we introduce the  definitions of the  Laplacian tensor   and the signless Laplacian tensor
   \textcolor{black}{for a $k$-uniform weighted hypergraph.}
 Let  $\mathbb{Q}_{k,n}$ be the set of
 {the connected}
 $k$-uniform weighted hypergraphs with $n$ vertices, where $k,n\geq 3$.
 Let  $G\in \mathbb{Q}_{k,n}$, where $k,n \geq 3$.
 We use $\mathcal{D}(G)=(d_{i_1i_2\cdots i_k})$ to denote a diagonal tensor of order $k$ and dimension $n$, where $k,n\geq 3$, $d_{i\ldots i}=w_{v_i}$ for $i\in [n]$ and $d_{i_1,\ldots i_k}=0$ otherwise.
 Let $\mathcal{L}(G) = \mathcal{D}(G) - \mathcal{A}(G)$ and
 $\mathcal{Q}(G) = \mathcal{D}(G) + \mathcal{A}(G)$.
 We call  $\mathcal{L}(G)$  and  $\mathcal{Q}(G)$ the  Laplacian tensor   and the signless Laplacian tensor of $G$, respectively.

   The research on the  spectral radius of the adjacency tensor, the Laplacian tensor and the signless Laplacian tensor for hypergraphs has attracted a lot of
   interests. For the three tensors of hypergraphs, many interesting results  about the characterization of the hypergraph with extremal spectral radius are derived,  and some properties and  bounds for the  extremal spectral radii have been obtained.
   Interested readers can find  Refs.
    \cite{2016-Yuan-p206,2020-Wen-HuanWang-FMC,--p,2017-PengXiao-p33,2018-Zhang-p1489,
  2018-HaiYanGuo-p236,2020-Wang-p,2019-Xiao-p1392,2016-YiZhengFan-p845,2018-LiYingKang-p661,2017-ChenOuyang-p141,Li-2013-p1001,2013-Xie-p1030,2014-Qi-p1045,2015-Hu-p1,2013-Hu-p2980,2016-Yue-p623,2020-Yue,2020-Yue-p2040007,2016-Yuan-p18}.

  Xie and Chang \cite{2013-Xie-p2195} obtained some bounds on the smallest and the largest $Z$-eigenvalues of the adjacency tensor for uniform hypergraphs.
    Xie and Chang \cite{2013-Xie-p1030,Xie-2013-p107} introduced  the signless Laplacian tensor for even uniform hypergraphs,
   and derived several properties of the smallest and the largest  $H$-eigenvalues and $Z $-eigenvalues of the signless Laplacian tensor for an even uniform hypergraph.
    Qi \cite{2014-Qi-p1045} defined the Laplacian and the signless Laplacian tensors of a uniform hypergraph for the study on  their $H^+$-eigenvalues and $H^{++}$-eigenvalues,
   and established some bounds for the largest signless Laplacian $H^+$-eigenvalue.
   Hu et al. \cite{2015-Hu-p1}   obtained a tight lower bound for the largest Laplacian $H$-eigenvalue of a $k$-uniform
  hypergraph and  derived the tight lower and upper bounds for the largest signless Laplacian $H$-eigenvalue   of a connected hypergraph.
     Yue et al. \cite{2016-Yue-p623}  obtained  the  upper bounds  of the largest Laplacian $H$-eigenvalue
     for a $k $-uniform loose path with a length not less than  3.
      All the results are related with the unweighted hypergraph.

 Inspired by the above results, in this paper, we investigate   the $H$-eigenvalue, the $H^{+}$-eigenvalue
  and the $H^{++}$-eigenvalue of adjacency tensor, Laplacian tensor and signless Laplacian tensor for the $k$-uniform weighted hypergraph $G$.
 This article is organized as follows.
 In Section 2,  some  notations and necessary lemmas which are useful for subsequent  proofs are introduced and  some basic properties of  the eigenvalues of  $\mathcal{A}(G)$, $\mathcal{L} (G)$, and $\mathcal{Q} (G)$ are presented.
 In Section 3, we study the  lower and upper bounds of the $H$-eigenvalue, the $H^{+}$-eigenvalue and the $H^{++}$-eigenvalue for $\mathcal{L} (G)$.
 The largest $H^+$-eigenvalue of $\mathcal{L} (G)$ is characterized.
 A relationship among the $H$-eigenvalues of $\mathcal{L} (G)$, $\mathcal{Q} (G)$ and $\mathcal{A} (G)$ is also given.
 In Section 4, we consider the lower and upper bounds of the $H$-eigenvalue, the $H^{+}$-eigenvalue and
 {the spectral radius}
 of $\mathcal{Q} (G)$.  The smallest $H^+$-eigenvalue of $\mathcal{Q} (G)$ is characterized.
 A property of the $H^{+}$-eigenvalue of $\mathcal{Q} (G)$ is derived.
 Finally, the lower and upper bounds of the $H$-eigenvalue, the $H^{+}$-eigenvalue and
 {the spectral radius}
 of $\mathcal{
  A} (G)$ are derived  in Section 5.
  A property of the $H^{+}$-eigenvalue for $\mathcal{A} (G)$ is also deduced.

\textcolor{black}{\section{Preliminary}}

   In this section, we first define some notations and introduce necessary lemmas. Then we
   derive some fundamental properties  about the eigenvalues of   $\mathcal{A}(G)$, $\mathcal{L} (G)$, and $\mathcal{Q} (G)$
    \textcolor{black}{for a $k$-uniform weighted hypergraph $G$.}

  Let  $\bm{x}=(x_{1},\cdots,x_{n})^{\textrm{T}}$   be an $n$-dimensional eigenvector of $\mathcal{A}(G)$ ($\mathcal{L} (G)$ and $\mathcal{Q} (G)$) and $x_i$ the component of $\bm{x}$ which corresponds to vertex $v_i$ $(i=1,\ldots,n)$ of $G$, where $G\in \mathbb{Q}_{k,n}$ with $k,n\geq 3$.
   Let  $U$ be a subset of $[n]$.
    Let  $$x^U = \prod_{i\in U} x_i.$$
 By (\ref{444444}) and (\ref{A}), for $i\in[n]$,  we get
 \begin{equation}\label{1.1}
 \left(\mathcal{A} (G) \bm{x}\right)_i = \sum\limits_{e\in E_{G}(v_i)} w_{G}(e) x^{e\setminus\{v_i\}}.
 \end{equation}
 Furthermore, we have
 \begin{equation}\label{1.2}
 \left(\mathcal{Q} (G) \bm{x}\right)_i = w_{v_i} x^{k-1}_i + \sum\limits_{e\in E_{G}(v_i)} w_{G}(e) x^{e\setminus\{v_i\}},
 \end{equation}
 \begin{equation}\label{1.3}
\left(\mathcal{L} (G) \bm{x}\right)_i = w_{v_i} x^{k-1}_i - \sum\limits_{e\in E_{G}(v_i)} w_{G}(e) x^{e\setminus\{v_i\}}.
\end{equation}

 Friedland et al. \cite{2013-Friedland-p738} defined the nonnegative weakly irreducible tensor
 and Yang et al. \cite{-N-p} restated it as follows.

\begin{definition}\label{definition2.3}{{\cite{-N-p} }}
 Let $\mathcal{A}=(a_{i_1i_2\cdots i_k})$ be a nonnegative tensor of order $k$ and dimension  $n$.
 If for any nonempty proper index subset $I \subset [n]$, there is at least an entry $a_{i_1i_2\cdots i_k}>0$,
 where $i_{1}\in I$ and at least an $i_{j}\in {{[n]\setminus I}}$ for $j=2,3, \ldots, k$,
 then $\mathcal{A}$ is called a nonnegative weakly irreducible tensor.
  \end{definition}

 \begin{lemma}\label{lemma2.2}\cite{2013-Friedland-p738,2010-Yang-p2517}
 Let $\mathcal{A}$ be a nonnegative tensor of order $k$ and dimension $n$, where $k\geq 2$.
 Then we have the following statements.

 (i). $\rho(\mathcal{A})$ is an eigenvalue of $\mathcal{A}$ with a nonnegative eigenvector $\bm{x}\in \mathbb{R}^{n}_{+}$ corresponding to it.

 (ii). If $\mathcal{A}$ is weakly irreducible,
 then $\rho(\mathcal{A})$ is the only eigenvalue of $\mathcal{A}$ with a positive eigenvector $\bm{x}\in \mathbb{R}^{n}_{++}$, up to a positive scaling coefficient.
 \end{lemma}

 Let  $\mathbb{S}_{k,n}$ be the set of real symmetric tensors of order $k$ and dimension $n$, where $k,n\geq3$.

 \begin{lemma}\label{lemma3.1}\cite{2005-Qi-p1302}
 We have the following conclusions on the eigenvalues of $\mathcal{A} \in \mathbb{S}_{k,n}$, {where $k,n\geq3$}.

 (i). A number $\lambda\in \mathbb{C}$ is an eigenvalue of $\mathcal{A}$ if and only if it is a root of the characteristic
 polynomial $\phi(\lambda)=\det(\mathcal{A}-\lambda\mathcal{I})$, where $\mathcal{I}$ is the unit tensor.

 (ii). The number of eigenvalues of $\mathcal{A}$ is $ n(k-1) ^{n-1} $. Their product is equal to
 $\det(\mathcal{A})$.

 (iii). The sum of all the eigenvalues of $\mathcal{A}$ is $(k - 1) ^{n-1} \textrm{tr}(\mathcal{A})$.

 \end{lemma}

 \begin{lemma}\label{lemma3.2}\cite{2005-Qi-p1302}
  Let $\mathcal{A}=(a_{i_{1}i_{2}\ldots i_{k}}) \in \mathbb{S}_{k,n}$, {where $k,n\geq3$}.  The following conclusions hold for $\mathcal{A}$.

  (i). Assume that $k$ is even.  $\mathcal{A}$ always has $H$-eigenvalues.  $\mathcal{A}$ is positive definite (positive semi-definite) if and only
  if all of its $H$-eigenvalues are positive (nonnegative).

 (ii).  The eigenvalues of $\mathcal{A}$ lie in the union of $n$ disks in $\mathbb{C}$.
  These $n$ disks have
 the diagonal elements $a_{i,\ldots,i}$ of the supersymmetric tensor as their centers, and the sums of the
 absolute values of the off-diagonal elements $\sum\limits_{i_2 ,\ldots, i_k=1; \{i_2 ,\ldots, i_k\}\neq\{i,\ldots,i\}}^{n} |a_{ii_{2}\ldots i_{k}}| $ as their radii, where $i\in[n]$.
 \end{lemma}

    Let $\mathcal{A}=(a_{i_{1}i_{2}\ldots i_{k}})\in \mathbb{R}^{n\times\cdots\times n}$ be a nonnegative tensor of  order $k$ and dimension $n$.
 Based on  $\mathcal{A}$, we define a
 \textcolor{black}{directed}
 graph  $\Gamma_{\mathcal{A}}$ as follows.
 The vertex set of   $\Gamma_{\mathcal{A}}$ is
 $V(\Gamma_{\mathcal{A}}) = \{1,\ldots,n\}$ and the arc set of $\Gamma_{\mathcal{A}}$ is
  \begin{equation}\label{PPP}
 E(\Gamma_{\mathcal{A}}) = \{(i,j) : a_{ii_{2}\ldots i_{k}}> 0, j\in\{i_{2}\ldots i_{k}\} \}.
 \end{equation}
  A graph is strongly connected if it contains a directed path from $i$ to $j$ and a directed path from $j$ to $i$ for every pair of vertices
 $i$ and $j$.  A tensor $\mathcal{A}$ is called weakly irreducible if $\Gamma_{\mathcal{A}}$ is strongly connected \cite{2013-Friedland-p738, 2014-Pearson-p1233, 2015-Bu-p168}.
 
  According to the definitions of the weakly irreducible tensor  and the adjacency tensor  of weighted hypergraph, we have Theorem \ref{theorem2.33} as follows.

 \begin{theorem}\label{theorem2.33}
 Let  $G\in \mathbb{Q}_{k,n}$, where $k,n \geq 3$.
  Any two of the three conclusions are equivalent.

 (i). $G$ is connected.

 (ii). $\mathcal{A}(G)$ is weakly irreducible.

 (iii). $\mathcal{Q}(G)$ is weakly irreducible.
 \end{theorem}

\noindent\textbf{Proof}.  Let  $G\in \mathbb{Q}_{k,n}$, where $k,n \geq 3$.
 Let the directed graph  associated with  $G$ be $\Gamma_{\mathcal{A}(G)} = \left(V(\Gamma_{\mathcal{A}} (G)), E(\Gamma_{\mathcal{A}} (G)) \right)$,
 where $V(\Gamma_{\mathcal{A}} (G))= \{ 1, 2,\ldots, n\}$ and
 $E(\Gamma_{\mathcal{A}} (G)) = \left\{(i,j) :
 e = \{v_i, v_j, v_{j_{3}},\ldots ,v_{j_{k}}\} \in E(G), j_3, \ldots, j_k\in[n]\backslash\{i,j\} \right\}$
 (by (\ref{A}) and (\ref{PPP})).

  Let $i$ and $j$ be any two different vertices in $V(\Gamma_{\mathcal{A}} (G))$.
  By the definition of $E(\Gamma_{\mathcal{A}} (G))$,
   for $i,j\in[n]$ and $i\neq j$, we obtain
 \begin{equation}\label{AAAA}
 (i,j),(j,i)\in E(\Gamma_{\mathcal{A}} (G))\Leftrightarrow v_i, v_j \in e \in E(G) .
 \end{equation}
  {Since $G$ is connected,} by (\ref{AAAA}), we can get that $\Gamma_{\mathcal{A}(G)}$ is  strongly connected.
  Namely, $\mathcal{A}(G)$ is weakly irreducible.

 If $\Gamma_{\mathcal{A}(G)}$ is  strongly connected, then for any two different vertices $i$ and $j$ in $V(\Gamma_{\mathcal{A}} (G))$, there exist  $j_1,\ldots,j_t\in V(\Gamma_{\mathcal{A}} (G))$ such that $(i,j_1),(j_1, j_2),\ldots,(j_t, j),(j, j_t),\ldots,(j_2, j_1),(j_1,\\
 i)\in E(\Gamma_{\mathcal{A}} (G))$, where $t\geq 0$.
 It follows from (\ref{AAAA}) that there exist   $e_1, e_2,\ldots,e_{t+1}\in E(G)$ such that $v_i, v_{j_1} \in e_1, \quad v_{j_1}, v_{j_2}\in e_2, \quad\ldots, \quad v_{j_t}, v_j\in e_{t+1}$, where $t\geq 0$.
 Namely, for $v_i, v_j\in V(G)$,  there exists a path in $G$ connecting
 $v_i$ and $v_j$. Thus, we get that $G$ is connected.

 Therefore, we have
 $(i)\Leftrightarrow(ii)$. Similarly, we can get $(i)\Leftrightarrow(iii)$. Thus, we have Theorem \ref{theorem2.33}.
 ~~$\Box$

 \begin{theorem}\label{theorem2.5}
 Let  $G\in \mathbb{Q}_{k,n}$, where $k,n \geq 3$.
 \textcolor{black}{$\rho(\mathcal{A}(G))$}
 %{Then $\rho(\mathcal{A(G)})$}
 $( \rho(\mathcal{Q}(G)) )$ is the only eigenvalue of $\mathcal{A}(G)$ $( \mathcal{Q}(G))$ with a unique positive eigenvector $\bm{x}\in \mathbb{R}^{n}_{++}$,
 up to a positive scaling coefficient.
 \end{theorem}

 \noindent\textbf{Proof}.  Let  $G\in \mathbb{Q}_{k,n}$, where $k,n \geq 3$.
  Since $G$ is connected, by  Theorem \ref{theorem2.33}, $\mathcal{A} (G)$ and
 $\mathcal{Q}(G)$ are weakly irreducible. Furthermore, by Lemma \ref{lemma2.2}(ii), we get Theorem \ref{theorem2.5}.
~~$\Box$

  By Lemmas \ref{lemma3.1} and \ref{lemma3.2}, we  obtain some basic properties of the eigenvalues of  $\mathcal{A} (G)$, $\mathcal{L}(G)$, and
 $\mathcal{Q}(G)$, where $G\in \mathbb{Q}_{k,n}$, which are shown in Theorem \ref{theorem3.1}.

 \begin{theorem}\label{theorem3.1}
 Let $G\in \mathbb{Q}_{k,n}$, where  $k,n \geq 3$.
 We have the five conclusions as follows.

(i). A number $\lambda\in \mathbb{C}$  is  an eigenvalue of $\mathcal{A}(G)$ ($\mathcal{L} (G)$ and $\mathcal{Q} (G)$) if and only if it is a root of the characteristic  polynomial $\phi(\mathcal{A}(G))$ ($\phi(\mathcal{L}(G))$ and $\phi(\mathcal{Q}(G))$).

 (ii).  The number of the eigenvalues of $\mathcal{A}(G)$ ($\mathcal{L} (G)$ and $\mathcal{Q} (G)$) is $ n(k-1) ^{n-1} $.
 Their product is equal to $\det(\mathcal{A}(G))$ ($\det(\mathcal{L} (G))$ and $\det(\mathcal{Q} (G))$).

(iii). The sum of all the eigenvalues of $\mathcal{A}(G)$  is zero and the sum of all the eigenvalues of  $\mathcal{L} (G)$ and $\mathcal{Q} (G)$ is  $(k-1)^{n-1} \sum\limits_{i=1}^{n} w_{v_i} = k(k-1)^{n-1} \sum\limits_{e\in E(G)} w_{G}(e)$.

(iv). All the eigenvalues of $\mathcal{A}(G)$ lie in the  disks $\{ \lambda :| \lambda | \leq W_0 \bigtriangleup\}$
 and all the eigenvalues of $\mathcal{L} (G)$ and $\mathcal{Q} (G)$ lie in the  disks $\{ \lambda :| \lambda - W_0 \bigtriangleup | \leq W_0 \bigtriangleup \}$.

(v).  When $k$ is even, $\mathcal{L} (G)$ and $\mathcal{Q} (G)$ are positive semi-definite tensors.
\end{theorem}

\noindent\textbf{Proof}. (i). The proof of Theorem \ref{theorem3.1}(i)--(iii).

  By Lemma \ref{lemma3.1}(i)--(iii) and the definition of the tensor of the weighted hypergraph, we can directly get  Theorem \ref{theorem3.1}(i)--(iii), respectively.

 (ii). The proof of Theorem \ref{theorem3.1}(iv).

  By (\ref{A}), for $i\in[n]$,  we have
  $$\sum\limits_{i_2 ,\ldots, i_k=1; \{i_2 ,\ldots, i_k\}\neq\{i,\ldots,i\}}^{n} a_{i i_{2}\ldots i_{k}} = \sum\limits_{e\in E_{G}(v_i)} w_{G}(e) = w_{v_i}.$$
 Let $\lambda$ be an arbitrary eigenvalue of $\mathcal{A}(G)$.
 Let $\bigcirc_i=\{\lambda :| \lambda | \leq  w_{v_i}\}$ be a disk, where $i=1,\ldots, n$.
 By  Lemma \ref{lemma3.2}(ii), we obtain $\lambda\in \bigcup_{i=1}^{n} \bigcirc_i$.
  Let $e$ be an arbitrary edge in $E(G)$.
  Since $w_{G}(e) \leq  W_0$ and $|E_{G}(v_i)| \leq \bigtriangleup$ for $i=1,\ldots, n$,  we get
  \textcolor{black}{\begin{align}
  |\lambda | \leq w_{v_i} = \sum\limits_{e\in E_{G}(v_i)} w_{G}(e)\leq W_0 \sum\limits_{e\in E_{G}(v_i)} 1
  = W_0 |E_{G}(v_i)|  \leq W_0\bigtriangleup.
       \end{align}}

 Let $\mu$ be an arbitrary eigenvalue of $\mathcal{L} (G)$ ($\mathcal{Q}(G)$).
 Let $\odot_i = \{\mu :| \mu - w_{v_i} | \leq w_{v_i} \}$ be a disk, where $i=1,\ldots, n$.
 By  Lemma \ref{lemma3.2}(ii), we get $\mu\in \bigcup_{i=1}^{n} \odot_i$.
 Since  $w_{v_i}  = \sum\limits_{e\in E_{G}(v_i)} w_{G}(e) \leq W_0\bigtriangleup$,
 we have $\mu \in \bigcup_{i=1}^{n} \odot_i\subseteq \{\mu :| \mu - W_0\bigtriangleup | \leq W_0\bigtriangleup\}$.

 (iii). The proof of Theorem \ref{theorem3.1}(v).

 When $k$ is even, by Theorem \ref{theorem3.1}(iv) and Lemma \ref{lemma3.2}(i), we obtain that $\mathcal{L} (G)$ and $\mathcal{Q} (G)$ are positive semi-definite tensors.
~~$\Box$

\section{The eigenvalues of $\mathcal{L} (G)$}

  In this section, we study the eigenvalues of $\mathcal{L} (G)$, where $G\in \mathbb{Q}_{k,n}$ {with $k,n\geq3$}.
  We obtain the upper and lower bounds of the $H$-eigenvalue and the $H^{+}$-eigenvalue of $\mathcal{L} (G)$, which are shown in  Theorems \ref{theorem3.3} and \ref{theorem3.7}, respectively.
  The  largest $H^+$-eigenvalue of $\mathcal{L} (G)$ is given in Theorem \ref{theorem4.2}.
  Two results of the $H^{+}$-eigenvalue and the $H^{++}$-eigenvalue of $\mathcal{L} (G)$ are derived in  Theorems \ref{theorem4.1}  and \ref{theorem4.3}, respectively.
   A relationship among the $H$-eigenvalues of $\mathcal{L} (G)$, $\mathcal{Q} (G)$ and $\mathcal{A} (G)$  is shown in
  Theorem \ref{theorem3.4}.

\begin{theorem}\label{theorem3.3} (\textbf{The bounds for the $H$-eigenvalue of $\mathcal{L} (G)$})
 Let  $G\in \mathbb{Q}_{k,n}$, where $k,n \geq 3$.
\textcolor{black}{Then $\mathcal{L} (G)$ has an $H$-eigenvalue}
  $\lambda$ and $ 0\leq \lambda\leq 2W_0\bigtriangleup$.
\end{theorem}

\noindent\textbf{Proof}. By (\ref{1.3}) and $\left(\mathcal{L} (G) \bm{x}\right)_i = \lambda x^{k-1}_i$ ($i=1,\ldots, n$), we get
\begin{align}\label{ww}
 \left(\mathcal{L} (G) \bm{1}\right)_i = w_{v_i} - \sum\limits_{e\in E_{G}(v_i)} w_{G}(e) = 0 =0\cdot[\bm{1}]_i.
 \end{align}
 By (\ref{ww}), zero is an $H^{++}$-eigenvalue of $\mathcal{L} (G)$ and $\bm{1}$ is the eigenvector of $\mathcal{L} (G)$ corresponding to zero.
 {Thus,  $\mathcal{L}(G)$ has $H$-eigenvalues.}
  Let $\lambda$ be an $H$-eigenvalue of $\mathcal{L} (G)$. By Theorem \ref{theorem3.1}(iv), $ 0\leq \lambda\leq 2W_0\bigtriangleup $.
 Therefore, we obtain Theorem \ref{theorem3.3}.
~~$\Box$

\begin{theorem}\label{theorem3.7} (\textbf{The bounds for the $H^+$-eigenvalue of $\mathcal{L} (G)$})
  Let  $G\in \mathbb{Q}_{k,n}$, where $k,n \geq 3$.
\textcolor{black}{Then $\mathcal{L} (G)$ has an $H^+$-eigenvalue}
 $\lambda$ and $0\leq \lambda \leq \alpha.$
\end{theorem}

\noindent\textbf{Proof}. Let  $G\in \mathbb{Q}_{k,n}$, where $k,n \geq 3$.
   Since  $\mathcal{L} (G)$ has an $H^{++}$-eigenvalue zero (by (\ref{ww})),
    $\mathcal{L} (G)$ has $H^+$-eigenvalues.
 Let $\lambda$ be an    $H^+$-eigenvalue of $\mathcal{L} (G)$.
 By Theorem \ref{theorem3.3}, $\lambda\geq 0$.
 Let $\bm{x}$ be an $H^+$-eigenvector  of $\mathcal{L} (G)$ corresponding to $\lambda$.
 Then $\bm{x}\in \mathbb{R}^{n}_{+}$.
 Thus, the largest  component of $\bm{x}$ is positive.
  Without loss of generality, we  assume the largest  component of $\bm{x}$ is 1.
 Let $u\in V(G)$ and $x_u = 1$.
 Therefore, by (\ref{1.3}) and $\left(\mathcal{L} (G) \bm{x}\right)_u = \lambda x^{k-1}_u$,
 we have
$$w_{u} - \lambda = \sum\limits_{e \in E_G(u)} w_{G}(e) x^{e\setminus\{u\} }. $$
 Since $\bm{x}\in \mathbb{R}^{n}_{+}$ and $ w_{G}(e)>0$ for any edge $e$ in $E(G)$,
 we have $w_{u} - \lambda \geq 0$. Namely, $\lambda \leq w_{u}$.
 Since $ w_u \leq \alpha$, we get $0\leq \lambda \leq \alpha.$
~~$\Box$

 Let $\bm{e}^{(i)}$ be an $n $-dimensional vector satisfying $e^{(i)}_j = 1$ if $j=i$ and $e^{(i)}_j = 0$ if
 $j\neq i$, where $i, j = 1, 2, \ldots, n$.

 \begin{theorem}\label{theorem4.1} (\textbf{$H^+$-eigenvalue of $\mathcal{L} (G)$})
  Let  $G\in \mathbb{Q}_{k,n}$, where $k,n \geq 3$.
  Then for any $i\in [n]$, $w_{v_i}$ is  an $H^+$-eigenvalue of $\mathcal{L} (G)$ and
  \textcolor{black}{$\bm{e}^{(i)}$ is an $H^+$-eigenvector}
   of $\mathcal{L} (G)$ corresponding to $w_{v_i}$.
 \end{theorem}

\noindent\textbf{Proof}.
  Let $i,j\in [n]$. If  $j=i$,  since  $e^{(i)}_j = 1$, we get
\begin{align}
 \left(\mathcal{L} (G) \bm{e}^{(i)}\right)_j
  %&= w_{v_i} (e^{(i)}_j)^{k-1} - \sum\limits_{e \in E_G(v_j)} w_{G}(e) (e^{(i)})^{e\setminus \{v_j\} }\nonumber\\
   &= \left( \mathcal{D} (G) \bm{e}^{(i)}\right)_j -  \left(\mathcal{A} (G) \bm{e}^{(i)}\right)_j \nonumber\\
  &= \sum\limits_{j_2 ,\ldots, j_k=1}^{n} d_{j j_2 \ldots j_{k}}e^{(i)}_{j_2} \cdots e^{(i)}_{j_k}- \sum\limits_{j_2 ,\ldots, j_k=1}^{n} a_{j j_2 \ldots j_{k}}e^{(i)}_{j_2} \cdots e^{(i)}_{j_k} \nonumber\\
  &= d_{j j \ldots j}e^{(i)}_j \cdots e^{(i)}_j \nonumber\\
  &= d_{i i \ldots i} = w_{v_i} = w_{v_i} e^{(i)}_j.  \nonumber
     \end{align}

  If  $j\neq i$, since $e^{(i)}_j = 0$, we obtain
 \begin{align}
 \left(\mathcal{L} (G) \bm{e}^{(i)}\right)_j
  %&= w_{v_i} (e^{(i)}_j)^{k-1} - \sum\limits_{e \in E_G(v_j)} w_{G}(e) (e^{(i)})^{e\setminus \{v_j\} }\nonumber\\
   &= \left( \mathcal{D} (G) \bm{e}^{(i)}\right)_j -  \left(\mathcal{A} (G) \bm{e}^{(i)}\right)_j \nonumber\\
  &= \sum\limits_{j_2 ,\ldots, j_k=1}^{n} d_{j j_2 \ldots j_{k}}e^{(i)}_{j_2} \cdots e^{(i)}_{j_k}- \sum\limits_{j_2 ,\ldots, j_k=1}^{n} a_{j j_2 \ldots j_{k}}e^{(i)}_{j_2} \cdots e^{(i)}_{j_k} \nonumber\\
  &= d_{j j \ldots j}e^{(i)}_j \cdots e^{(i)}_j  \nonumber\\
  &= 0 = w_{v_i} e^{(i)}_j.  \nonumber
     \end{align}
 Therefore, we get Theorem \ref{theorem4.1}.
~~$\Box$

\begin{theorem}\label{theorem4.2}  (\textbf{The largest $H^+$-eigenvalue of $\mathcal{L} (G)$})
  Let  $G\in \mathbb{Q}_{k,n}$, where $k,n \geq 3$.
 Then $\alpha$ is the largest $H^+$-eigenvalue of $\mathcal{L} (G)$.
\end{theorem}

\noindent\textbf{Proof}. By  Theorem \ref{theorem3.7}, $\mathcal{L} (G)$ has $H^+$-eigenvalues.
 Let $\lambda$ be an $H^+$-eigenvalue of $\mathcal{L} (G)$.
 It  follows from Theorem \ref{theorem3.7} that $\lambda \leq \alpha$.
 Let $\lambda_{0} $ be the largest  $H^+$-eigenvalue of $\mathcal{L} (G)$.
 Then, we have $\lambda_{0} \leq \alpha$.
 By Theorem \ref{theorem4.1},
 $\alpha$ is an $H^+$-eigenvalue of $\mathcal{L} (G)$. Therefore, $\alpha\leq\lambda_{0}$.
 Thus, we obtain $\lambda_{0}=\alpha$.
~~$\Box$

\begin{theorem}\label{theorem4.3} (\textbf{$H^{++}$-eigenvalue of $\mathcal{L} (G)$})
  Let  $G\in \mathbb{Q}_{k,n}$, where $k,n \geq 3$.
 Then zero is the unique $H^{++}$-eigenvalue of $\mathcal{L} (G)$.
\end{theorem}

\noindent\textbf{Proof}.  By  (\ref{ww}), zero is an  $H^{++}$-eigenvalue of $\mathcal{L} (G)$.
 Thus, $\mathcal{L} (G)$ has $H^{++}$-eigenvalues.
 Let $\lambda$ be an $H^{++}$-eigenvalue of $\mathcal{L} (G)$.
 Next, we prove $\lambda=0$.

 Let $\bm{x}$ be an $H^{++}$-eigenvector of $\mathcal{L} (G)$ corresponding to $\lambda$.
  Then $\bm{x}\in \mathbb{R}^{n}_{++}$.
 Without loss of generality, we assume that the smallest component of $\bm{x}$ is 1.
 Let $v\in V(G)$ and $x_v = 1$.
 Therefore, by (\ref{1.3}) and $\left(\mathcal{L} (G) \bm{x}\right)_v = \lambda x^{k-1}_v$, we have
\begin{align}
 \lambda=  w_{v} - \sum\limits_{e \in E_G(v)} w_{G}(e) x^{e\setminus\{v\} }
 \leq w_{v} - \sum\limits_{e \in E_G(v)} w_{G}(e) \label{6.2}
 = 0.
\end{align}
 It is noted that (\ref{6.2}) holds since  $x_{v'} \geq 1$ for any $ v'\in V(G)$ and $ w_{G}(e)>0$ for any edge $e\in E(G)$.
 Therefore, we obtain $\lambda\leq0$.
 Furthermore, by Theorem \ref{theorem3.3},  $\lambda\geq0$. Therefore, we get  $\lambda=0$. Namely, zero is the unique $H^{++}$-eigenvalue of $\mathcal{L} (G)$.
 ~~$\Box$

Let $G^{\vartriangle}$ be a $k$-uniform hypergraph.
 It is interesting that Qi \cite{2014-Qi-p1045} used the methods of
 optimization theory to obtain a result  about the largest  $H^+$-eigenvalue
 of $\mathcal{L}(G^{\vartriangle})$ (namely,
 \textcolor{black}{Theorem 5.1}
  in  \cite{2014-Qi-p1045}) which is similar  to  Theorem \ref{theorem4.2},  and
   obtained some results  about the $H^+$-eigenvalue
 and the unique $H^{++}$-eigenvalue of $\mathcal{L}(G^{\vartriangle})$ (namely,
 \textcolor{black}{Theorem 3.2}
 in \cite{2014-Qi-p1045}) which are similar  to Theorems \ref{theorem4.1} and \ref{theorem4.3}.

  For $G\in \mathbb{Q}_{k,n}$ {with $k,n\geq3$}, we obtain a relationship of the $H$-eigenvalues of $\mathcal{A} (G)$, $\mathcal{L} (G)$ and $\mathcal{Q} (G)$, which is shown in Theorem \ref{theorem3.4}. We can prove Theorem \ref{theorem3.4} by using
  the relationship among the tensors of $\mathcal{A} (G)$, $\mathcal{L} (G)$ and $\mathcal{Q} (G)$.
  However, to enrich the diversity of proof, we prove it by using  different methods.

\begin{theorem}\label{theorem3.4} (\textbf{The relationship among the $H$-eigenvalues of $\mathcal{A} (G)$, $\mathcal{L} (G)$ and $\mathcal{Q} (G)$})
 Let  $G\in \mathbb{Q}_{k,n}$, where $k,n \geq 3$. Furthermore, we suppose that $G$ is $r$-regular ($r\geq1$) and each edge of $G$ has the same weight $W_0$.
 If $\lambda$ is an  $H$-eigenvalue of $\mathcal{L} (G)$, then

 (i). $ 2W_0 r -\lambda$ is an $H$-eigenvalue of $\mathcal{Q} (G)$.

 (ii). $ W_0 r -\lambda$ is an $H$-eigenvalue of $\mathcal{A} (G)$.
\end{theorem}

\noindent\textbf{Proof}.
 Let $G$ be as described in Theorem \ref{theorem3.4}.
 By Theorem \ref{theorem3.3}, $\mathcal{L} (G)$ has $H$-eigenvalues.
 Let $\lambda$ be an $H$-eigenvalue of $\mathcal{L} (G)$ and
 $\bm{x}$ be an $H$-eigenvector of $\mathcal{L} (G)$ corresponding to $\lambda$.
 Thus, we have $\bm{x}\in \mathbb{R}^{n}$.
 For any $i\in[n]$,  by (\ref{1.3}) and $\left(\mathcal{L} (G) \bm{x}\right)_i = \lambda x^{k-1}_i$, we get
\begin{align}
 \lambda x^{k-1}_i &= w_{v_i} x^{k-1}_i - \sum\limits_{e \in E_G(v_i)} w_{G}(e) x^{e\setminus \{v_i\} } \nonumber\\
 &= W_0\cdot r\cdot x^{k-1}_i - W_0\sum\limits_{e \in E_G(v_i)} x^{e\setminus \{v_i\} }.      \label{i}
\end{align}
 It it noted that (\ref{i}) holds since $w_{v_i} = \sum\limits_{e \in E_G(v_i)} w_{G}(e)$,
 $G$ is $r$-regular,  and each edge of $G$ has weight $W_0$.

\textcolor{black}{For any $i\in[n]$, we have
\begin{align}
 (2W_0 r -\lambda) x^{k-1}_i &= 2 W_0\cdot r\cdot x^{k-1}_i - \lambda x^{k-1}_i \nonumber\\
% &= 2 W_0\cdot r\cdot x^{k-1}_i - W_0\cdot r\cdot x^{k-1}_i + W_0\sum\limits_{e \in E_G(v_i)} x^{e\setminus \{v_i\} }\nonumber\\
 &= W_0\cdot r\cdot x^{k-1}_i + W_0\sum\limits_{e \in E_G(v_i)} x^{e\setminus \{v_i\} }     \label{iv}\\
 &= w_{v_i} x^{k-1}_i + \sum\limits_{e \in E_G(v_i)} w_{G}(e) x^{e\setminus \{v_i\} } \label{v}\\
 &= \left(\mathcal{Q} (G) \bm{x}\right)_i.      \nonumber
\end{align}
  It is noted that  (\ref{iv}) follows from (\ref{i}), and
(\ref{v}) holds since $w_{v_i} = \sum\limits_{e \in E_G(v_i)} w_{G}(e)$,
 $G$ is $r$-regular and each edge of $G$ has the same weight $W_0$.}

 Therefore, we obtain that  $ 2W_0 r -\lambda$ is an $H$-eigenvalue of $\mathcal{Q} (G)$. Namely, we get  Theorem \ref{theorem3.4}(i).
 By the methods similar to those for theorem \ref{theorem3.4}(i), we get  Theorem \ref{theorem3.4}(ii).
~~$\Box$

\section{The eigenvalues of $\mathcal{Q} (G)$}

 In this section, we investigate the eigenvalues of $\mathcal{Q} (G)$, where $G\in \mathbb{Q}_{k,n}$ with $k,n\geq 3$.
    The upper and lower bounds for the $H$-eigenvalue, the $H^+$-eigenvalue and { the spectral radius} of $\mathcal{Q} (G)$  are shown in Theorems \ref{theorem3.2}--\ref{theorem3.8}, respectively. The weighted hypergraph with the largest $H$-eigenvalue of $\mathcal{Q} (G)$ is also characterized in Theorem \ref{theorem3.2}.
    {A property of the $H^+$-eigenvalue of $\mathcal{Q} (G)$ is given in Theorem \ref{theorem4.1jia}.
    The smallest  $H^+$-eigenvalue of $\mathcal{Q} (G)$ is obtained in  Theorem \ref{theorem4.4}.}

\begin{theorem}\label{theorem3.2} (\textbf{The bounds for the $H$-eigenvalue of $\mathcal{Q} (G)$})
  Let  $G\in \mathbb{Q}_{k,n}$, where $k,n \geq 3$.
 \textcolor{black}{Then (i)  $\mathcal{Q} (G)$ has an $H$-eigenvalue}
 $\lambda$ and
 $0\leq\lambda\leq 2W_0\bigtriangleup$;
 (ii)  $2W_0\bigtriangleup$ is an $H$-eigenvalue of $\mathcal{Q} (G)$ if and only if $G$ is $\bigtriangleup$-regular and each edge of $G$ has the same weight $W_0$.
 \end{theorem}

\noindent\textbf{Proof}.
 Let  $G\in \mathbb{Q}_{k,n}$, where $k,n \geq 3$.
 Since $G$ is connected, by Theorem \ref{theorem2.5},
 $\rho(\mathcal{Q} (G))$ is an $H^{++}$-eigenvalue of $\mathcal{Q} (G)$.
 Thus, $\mathcal{Q} (G)$ has $H$-eigenvalues.
 Let $\lambda$ be an $H$-eigenvalue of $\mathcal{Q} (G)$. By Theorem \ref{theorem3.1}(iv), $0\leq\lambda\leq 2W_0\bigtriangleup$. Thus, we get Theorem \ref{theorem4.1jia}(i).

 Next, we prove  Theorem \ref{theorem3.2}(ii).

   If $G$ is $\bigtriangleup$-regular and each edge of $G$ has the same weight $W_0$, by (\ref{1.2}), {for $i\in [n]$,} we get
\begin{align}
  \left(\mathcal{Q} (G) \bm{1}\right)_i &=  w_{v_i}+ \sum\limits_{e\in E_{G}(v_i)} w_{G}(e)
  = 2\sum\limits_{e\in E_{G}(v_i)} w_{G}(e)
  = 2W_0\bigtriangleup
  = 2W_0\bigtriangleup\cdot[\bm{1}]_i.  \nonumber
     \end{align}
 Thus, $2W_0\bigtriangleup$ is an $H$-eigenvalue of $\mathcal{Q} (G)$ and $\bm{1}$ is the eigenvector of $\mathcal{Q} (G)$ corresponding to $ 2W_0\bigtriangleup$.

 We assume that  $2W_0\bigtriangleup$ is an $H$-eigenvalue of $\mathcal{Q} (G)$. Next, we prove that $G$ is $\bigtriangleup$-regular and each edge of $G$ has the same weight $W_0$.
 Let $\bm{x}=(x_{1},\cdots,x_{n})^{\textrm{T}} \in \mathbb{R}^{n}$ be the eigenvector of $\mathcal{Q} (G)$ corresponding to $ 2W_0\bigtriangleup$ with$\sum\limits_{i=1}^{n} x^{k}_i = 1$. Without loss of generality, let $| x_j | = \max\limits_{1\leq i \leq n} | x_i |$, where $j\in [n]$. Obviously, $| x_j |>0$. By (\ref{1.2}) and $\left(\mathcal{Q} (G) \bm{x}\right)_j = 2W_0\bigtriangleup x^{k-1}_j$, we have
\begin{align}
w_{v_j} x^{k-1}_j + \sum\limits_{e= \{v_j, v_{j_2} ,\ldots, v_{j_k}\} \in E(G)} w_{G}(e) x_{j_2} \cdots x_{j_k}= 2W_0\bigtriangleup x^{k-1}_j. \label{B}
\end{align}
 Since $w_{v_j} = \sum\limits_{e\in E_{G}(v_j)} w_{G}(e)$, we have  $w_{v_j} < 2W_0\bigtriangleup$.
 Thus, we get
\begin{align}
  2W_0\bigtriangleup - w_{v_j} &= |\sum\limits_{e= \{v_j, v_{j_2} ,\ldots, v_{j_k}\} \in E(G)} w_{G}(e) \frac{x_{j_2}}{x_j} \cdots \frac{x_{j_k}}{x_j} | \label{4}\\
  &\leq \sum\limits_{e= \{v_j, v_{j_2} ,\ldots, v_{j_k}\} \in E(G)} w_{G}(e) | \frac{x_{j_2}}{x_j}| \cdots | \frac{x_{j_k}}{x_j} | \label{5}\\
  &\leq  \sum\limits_{e \in E_{G} (v_j)} w_{G}(e)= w_{v_j}. \label{6}
       \end{align}
    It is noted that (\ref{4}) is obtained from (\ref{B}) by first subtracting $w_{v_j} x^{k-1}_j$ from both sides of (\ref{B}), then dividing $x^{k-1}_j$ at the same time,  and finally taking the modulus on both sides to get (\ref{4}).
       (\ref{5}) follows from the property of absolute value inequality and   (\ref{6}) follows from $| x_j | = \max\limits_{1\leq i \leq n} | x_i |$.
 Thus, we get $W_0\bigtriangleup \leq  w_{v_j} = \sum\limits_{e \in E_{G} (v_j)} w_{G}(e) $.

  For a fixed $j\in [n]$ with $| x_j |=\max\limits_{1\leq i \leq n} | x_i |$,
  if $| E_{G} (v_j) | < \bigtriangleup$ holds  or there exists one $e \in E_{G} (v_j)$ such that $ w_{G}(e) < W_0$,
  then $w_{v_j} < W_0\bigtriangleup$.
  Obviously, this  contradicts  $W_0\bigtriangleup \leq  w_{v_j}$.
  Therefore,  we have  $| E_{G} (v_j) | = \bigtriangleup$  for a fixed $j\in [n]$ with $| x_j |=\max\limits_{1\leq i \leq n} | x_i |$ and
   $w_{G}(e) = W_0$ for any  $e \in E_{G} (v_j)$.
 Thus, the two equalities in (\ref{5}) and (\ref{6}) hold simultaneously.
 Namely,  $ | x_v | = | x_j |=\max\limits_{1\leq i \leq n} | x_i |$ for any $v\in N_G(v_j)$, and
  $x^{e_1\setminus \{v_j\} }$ and $x^{e_2\setminus \{v_j\} }$ have the same symbol,
   where $e_1,e_2\in E_G(v_j)$. %$e_1, e_2 \in E_G(v_j)$, and
  Since $G$ is connected, there exists a path between every pair of vertices in $V(G)$. By repeatedly using the same analysis as above, we get that $w_{G}(e) = W_0$ for any $e \in E(G)$ and $| E_{G} (v_i) | = \bigtriangleup$ for $i = 1,\ldots,n$.
  Therefore, if {$G$ is connected and} $2W_0\bigtriangleup$ is an $H$-eigenvalue of $\mathcal{Q} (G)$, then $G$ is $\bigtriangleup$-regular and each edge of $G$ has the same weight $W_0$. Namely, Theorem \ref{theorem3.2}(ii) holds.
  ~~$\Box$

\begin{theorem}\label{theorem3.6} (\textbf{The bounds for the $H^+$-eigenvalue of $\mathcal{Q} (G)$})
  Let  $G\in \mathbb{Q}_{k,n}$, where $k,n \geq 3$.
\textcolor{black}{Then $\mathcal{Q} (G)$ has an $H^+$-eigenvalue}
  $\lambda$ and $\delta\leq \lambda \leq 2\alpha.$
\end{theorem}

\noindent\textbf{Proof}.
 Let  $G\in \mathbb{Q}_{k,n}$, where $k,n \geq 3$.
 {By Theorem \ref{theorem2.5},
 $\rho(\mathcal{Q} (G))$ is an $H^{++}$-eigenvalue of $\mathcal{Q} (G)$.}
 Thus, $\mathcal{Q} (G)$ has $H^+$-eigenvalues.
 Let $\lambda$ be an $H^+$-eigenvalue of $\mathcal{Q} (G)$. Let $\bm{x}$ be an $H^+$-eigenvector of $\mathcal{Q} (G)$ corresponding to $\lambda$. Then, $\bm{x}\in \mathbb{R}^{n}_{+}$. Thus,  the largest component of  $\bm{x}$ is positive.
 Without loss of generality, we assume that the largest component of $\bm{x}$ is 1. Let $u\in V(G)$ and $x_u = 1$. By (\ref{1.2}) and $\left(\mathcal{Q} (G) \bm{x}\right)_u = \lambda x^{k-1}_u$, we get
\begin{align}
w_{u} + \sum\limits_{e \in E_G(u)} w_{G}(e) x^{e\setminus\{u\} }= \lambda. \label{liu}
\end{align}
  Since $\bm{x}\in \mathbb{R}^{n}_{+}$,  $w_{G}(e) >0$ for any $e\in E(G)$,  and $0\leq x_v \leq1$ for any $v\in V(G)$,
 we obtain
\begin{align}
 0\leq \lambda- w_{u} \label{qi}
 = \sum\limits_{e \in E_G(u)} w_{G}(e) x^{e\setminus\{u\} } \leq \sum\limits_{e \in E_G(u)} w_{G}(e)
 = w_{u}.
\end{align}
 Thus, we get $w_{u}\leq \lambda \leq 2w_{u}$. Since $\delta\leq w_u \leq \alpha$, we obtain $\delta\leq \lambda \leq 2\alpha.$
~~$\Box$

\begin{theorem}\label{theorem3.8} (\textbf{The bounds for {the spectral radius} of $\mathcal{Q} (G)$})
  Let  $G\in \mathbb{Q}_{k,n}$, where $k,n \geq 3$.
{Then $2 \delta\leq \rho(\mathcal{Q} (G)) \leq 2\alpha.$}
\end{theorem}

\noindent\textbf{Proof}. Let  $G\in \mathbb{Q}_{k,n}$, where $k,n \geq 3$.
 {Since $G$ is connected, by Theorem \ref{theorem2.5}, $\rho(\mathcal{Q}(G))$ is the only $H^{++}$-eigenvalue of $\mathcal{Q} (G)$.}
  Let $\bm{x}$ be the $H^{++}$-eigenvector of  $\mathcal{Q} (G)$ corresponding to {$\rho(\mathcal{Q}(G))$}.
 We have $\bm{x}\in \mathbb{R}^{n}_{++}$. Let $v\in V(G)$ and $x_v$ be the smallest component of $\bm{x}$. By
 (\ref{1.2}) and
 {$\left(\mathcal{Q} (G) \bm{x}\right)_v = \rho(\mathcal{Q}(G)) x^{k-1}_v$,} we get
$$w_{v}x^{k-1}_v + \sum\limits_{e \in E_G(v)} w_{G}(e) x^{e\setminus\{v\} }= \rho(\mathcal{Q}(G)) x^{k-1}_v.$$
 Since $\bm{x}\in \mathbb{R}^{n}_{++}$, we have $x_v>0$.  Thus, we obtain
$$\rho(\mathcal{Q}(G))  = w_v + \sum\limits_{e \in E_G(v)} w_{G}(e) \frac{x^{e\setminus\{v\} }}{x^{k-1}_v}.$$
 Since $0< x_{v} \leq x_{v'}$ for any $v'\in V(G)$  and $ w_{G}(e)>0$ for any  $e\in E(G)$, we obtain $\rho(\mathcal{Q}(G)) \geq w_v + \sum\limits_{e \in E_G(v)} w_{G}(e)\geq 2 \delta$.
 {Furthermore, by Theorem \ref{theorem3.6}, we get $\rho(\mathcal{Q} (G)) \leq 2\alpha$.
 Thus, we obtain Theorem \ref{theorem3.8}.}
~~$\Box$

  By using the methods similar to those for Theorem \ref{theorem4.1}, we get Theorem \ref{theorem4.1jia}.

  \begin{theorem}\label{theorem4.1jia} (\textbf{$H^+$-eigenvalue of $\mathcal{Q} (G)$})
  Let  $G\in \mathbb{Q}_{k,n}$, where $k,n \geq 3$.
\textcolor{black}{Then for any $i\in [n]$, $w_{v_i}$ is an $H^+$-eigenvalue of $\mathcal{Q} (G)$
 and $\bm{e}^{(i)}$ is an $H^+$-eigenvector of  $\mathcal{Q} (G)$ corresponding to   $w_{v_i}$.}
 \end{theorem}

In Theorem \ref{theorem4.4}, we obtain the smallest  $H^+$-eigenvalue of $\mathcal{Q} (G)$.
 The proof of Theorem \ref{theorem4.4} is omitted since we can apply
 Theorems \ref{theorem3.6} and \ref{theorem4.1jia} and use the same methods similar to those for
 Theorem \ref{theorem4.2} to get it.

 \begin{theorem}\label{theorem4.4}  (\textbf{The smallest  $H^+$-eigenvalue of $\mathcal{Q} (G)$})
  Let  $G\in \mathbb{Q}_{k,n}$, where $k,n \geq 3$.
 Then $\delta$ is the smallest  $H^+$-eigenvalue of $\mathcal{Q} (G)$.
 \end{theorem}

 It is interesting that  Qi \cite{2014-Qi-p1045} used the methods of optimization theory to obtain a result  about the smallest  $H^+$-eigenvalue of the signless Laplacian tensor of a $k$-uniform hypergraph (shown in
 \textcolor{black}{Theorem 7.1}
 in  \cite{2014-Qi-p1045}) which is similar  to  Theorem \ref{theorem4.4}.

 \section{The eigenvalues of $\mathcal{A} (G)$}

 In this section, we study the eigenvalues of $\mathcal{A} (G)$, where $G\in \mathbb{Q}_{k,n}$ {with $k,n\geq3$}.
 The upper and lower bounds for the $H$-eigenvalue, the $H^{+}$-eigenvalue and {the spectral radius} of $\mathcal{A} (G)$ are  derived in Theorems \ref{theorem3.10}--\ref{theorem3.9}, respectively. The weighted hypergraph with the largest $H$-eigenvalue of $\mathcal{A} (G)$ is also presented  in Theorem \ref{theorem3.10}.
  We find that zero is an $H^+$-eigenvalue of $\mathcal{A} (G)$ and
  \textcolor{black}{$\bm{e}^{(i)}$ is an $H^+$-eigenvector }
    of $\mathcal{A} (G)$ corresponding to zero, where $i\in[n]$, which is shown in Theorem \ref{theorem4.5}.

\begin{theorem}\label{theorem3.10} (\textbf{The bound for the $H$-eigenvalue of $\mathcal{A} (G)$})
  Let  $G\in \mathbb{Q}_{k,n}$, where $k,n \geq 3$.
        \textcolor{black}{Then (i)  $\mathcal{A} (G)$ has an $H$-eigenvalue}
   $\lambda$ and  $|\lambda|\leq W_0\bigtriangleup$;
    (ii)   \textcolor{black}{$W_0\bigtriangleup$ is an $H$-eigenvalue}
    of $\mathcal{A} (G)$ if and only if $G$ is $\bigtriangleup$-regular and each edge of $G$ has the same weight $W_0$.
    \end{theorem}

 \noindent\textbf{Proof}. Let  $G\in \mathbb{Q}_{k,n}$, where $k,n \geq 3$. Since $G$  is connected, by Theorem \ref{theorem2.5},
 $\mathcal{A} (G)$ has  an  $H^{++}$-eigenvalue $\rho(G)$.
  Thus, $\mathcal{A} (G)$ has an $H$-eigenvalue $\lambda$.
 By Theorem \ref{theorem3.1}(iv), $ |\lambda| \leq W_0\bigtriangleup $. Therefore, we get Theorem \ref{theorem3.10}(i).
 By the methods similar to those for Theorem \ref{theorem3.2}, we obtain Theorem \ref{theorem3.10}(ii).
 ~~$\Box$

  Let $G\in \mathbb{Q}_{k,n}$, where $k,n \geq 3$.
  Let $X$ be a non-empty subset of $V(G)$.
  We use $E_t (X)$  to denote the set of edges of $G$ which share  $t$ common vertices with $X$, where $t\geq 1$.
 Namely, $E_t (X) = \{e: e\in E(G) \textrm{ and } |e\cap X |=t  \}$.
 Furthermore, let $E^{v}_t (X) = \{e: e\in E(G), v\in e ~\textrm{and}~ |e\cap X |=t  \}$.
 We define $e_G(u,v)$ as the number of the edges of $G$ which contain $u$ and $v$, where $u,v\in V(G)$.

 \begin{theorem}\label{theorem3.5} (\textbf{The bounds for the $H^+$-eigenvalue of $\mathcal{A} (G)$})
   Let  $G\in \mathbb{Q}_{k,n}$, where $k,n \geq 3$.
\textcolor{black}{Then $\mathcal{A} (G)$ has an $H^+$-eigenvalue}
   $\lambda$ and
  $$0\leq\lambda \leq \sqrt{\frac{W^2_0}{k-1} \sum^{k}\limits_{t=1}  \sum\limits_{e\in E_t (N_G(u))} \sum\limits_{v\in (e\cap N_G(u))} e_G(u,v)},
 $$
 where $u$ is the vertex of $G$ which has the largest component of the principal eigenvector corresponding to $\lambda$.
 \end{theorem}

 \noindent\textbf{Proof}.
  Let $G\in \mathbb{Q}_{k,n}$, where $k,n \geq 3$. % By  (\ref{A}) and Lemma \ref{lemma2.2}(i),
  {By Theorem \ref{theorem2.5},
 $\rho(G)$ is an $H^{++}$-eigenvalue of $\mathcal{A} (G)$.}
  Thus, $\mathcal{A} (G)$ has $H^+$-eigenvalues.
   Let $\lambda$ be an $H^+$-eigenvalue of $\mathcal{A} (G)$ and
  \textcolor{black}{$\bm{x}$ be an $H^+$-eigenvector}
   of $\mathcal{A} (G)$ corresponding to $\lambda$.
  Thus, $\bm{x}\in \mathbb{R}^{n}_{+}$.
  For all $i\in[n]$, by (\ref{1.1}) and $\left(\mathcal{A} (G) \bm{x}\right)_i = \lambda x^{k-1}_i$, we get
 \begin{align}\label{vii}
 \lambda x^{k-1}_i = \sum\limits_{e=\{v_{i},v_{i_2},\ldots,v_{i_k} \} \in E(G)} w_{G}(e) x_{i_2}\cdots x_{i_k}.
 \end{align}
 Since $\bm{x}\in \mathbb{R}^{n}_{+}$, the largest component of $\bm{x}$ is positive.
 Without loss of generality, we assume the largest component of $\bm{x}$ is 1.
  Let $u\in V(G)$ and $x_u = 1$.
 In (\ref{vii}), let $i=u$. Since $\bm{x}\in \mathbb{R}^{n}_{+}$ and $ w_{G}(e)>0$ for any $e\in E(G)$, by  (\ref{vii}), we obtain
 \begin{align}
 0\leq \lambda &= \sum\limits_{e=\{u,v_{i_2},\ldots,v_{i_k} \} \in E(G)} w_{G}(e) x_{i_2}\cdots x_{i_k} \label{6.4}\\
 &\leq \frac{1}{k-1}\sum\limits_{e=\{u,v_{i_2},\ldots,v_{i_k} \} \in E(G)} w_{G}(e) \left( x^{k-1}_{i_2} + \cdots + x^{k-1}_{i_k}\right), \label{viii}
 \end{align}
 where (\ref{viii}) follows from the AM-GM equality.

  Multiplying both sides of (\ref{viii}) by $\lambda $ and bearing (\ref{vii}) in mind, we have
 \begin{align}
 &\lambda^2  \leq \frac{1}{k-1}\sum\limits_{e=\{u,v_{i_2},\ldots,v_{i_k} \} \in E(G)} w_{G}(e)\quad \times \nonumber\\
 &\qquad\left[\sum\limits_{f_2 \in E_G(v_{i_2})} w_{G}(f_2) x^{f_2 \setminus\{v_{i_2} \}} + \cdots
  + \sum\limits_{f_k \in E_G(v_{i_k})} w_{G}(f_k) x^{f_k \setminus\{v_{i_k} \}}\right].\label{1}
 \end{align}
 Since $0\leq x_v \leq1$ for any $v\in V(G)$ and $w_{G}(e)\geq 0$ for any $e\in E(G)$, by (\ref{1}),  we get
 \begin{align}
 &\lambda^2\leq \frac{1}{k-1}\sum\limits_{e=\{u,v_{i_2},\ldots,v_{i_k} \} \in E(G)} w_{G}(e)\nonumber
 \left[\sum\limits_{f_2 \in E_G(v_{i_2})} w_{G}(f_2) + \cdots + \sum\limits_{f_k \in E_G(v_{i_k})} w_{G}(f_k)\right]\nonumber.
  \end{align}
 For any $f_{s}\in E_G(v_{i_s})$, since $v_{i_s} \in N_G(u)$,
 we have $1\leq|f_s\cap N_G(u)|   \leq k$, where $s=2,\ldots,k$.
  Thus, $E_G(v_{i_s}) = \bigcup\limits_{t=1} ^k E_{t}^{v_{i_s} } (N_G(u))$. Therefore, we obtain
  \begin{align}
 &\lambda^2\leq \frac{1}{k-1}\sum\limits_{e=\{u,v_{i_2},\ldots,v_{i_k} \} \in E(G)} w_{G}(e)\quad \times  \nonumber\\
 &\qquad \left[\sum\limits_{t=1}^{k}\left(\sum\limits_{f_2 \in E^{v_{i_2}}_t (N_G(u))} w_{G}(f_2)  + \cdots + \sum\limits_{f_k \in E^{v_{i_k}}_t (N_G(u))} w_{G}(f_k)\right)\right].\label{san}
 \end{align}
 Since $0< w_G(e)\leq W_0$ for any $e\in E(G)$, we have
  \begin{align}
 &\lambda^2 \leq \frac{W^2_0}{k-1}\sum\limits_{e=\{u,v_{i_2},\ldots,v_{i_k} \} \in E(G)}
 \left[\sum\limits_{t=1}^{k} \left(\sum\limits_{f_2 \in E^{v_{i_2}}_t (N_G(u))} 1  + \cdots
 + \sum\limits_{f_k \in E^{v_{i_k}}_t (N_G(u))} 1\right)\right]\nonumber.
 \end{align}
  Obviously, the upper bound of  $\lambda^2$ is related with these edges which contain at least one vertex in $N_G(u)$.
  Let $f=\{v_{j_1}, v_{j_2}, \ldots, v_{j_k} \} \in E(G)$ with $f\cap N_G(u) \neq\emptyset$.
   Without loss of generality, we suppose $f\cap N_G(u) = \{v_{j_1}, v_{j_2}, \ldots, v_{j_t} \}$, where $1\leq t\leq k$.
  Then $f$ appears $e_G(u,v_{j_1})+ e_G(u,v_{j_2})+ \cdots + e_G(u,v_{j_t})$ times. Therefore, we obtain
  \begin{align}
 &\lambda^2 \leq \frac{W^2_0}{k-1}\left[\sum\limits_{e\in E_1(N_G(u))} \sum\limits_{v \in (e\cap N_G(u))} e_G(u,v)+\cdots+
 \sum\limits_{e \in E_k (N_G(u))} \sum\limits_{v \in (e\cap N_G(u))} e_G(u,v)\right] \label{wu}\\
 &= \frac{W^2_0}{k-1}\sum\limits_{t=1}^{k} \sum\limits_{e \in E_t (N_G(u))} \sum\limits_{v \in (e\cap N_G(u))} e_G(u,v).\nonumber
\end{align}
 Thus, we get Theorem \ref{theorem3.5}. ~~$\Box$

 It is noted that Theorem \ref{theorem3.5} is
 \textcolor{black}{a generalization}
 of Lemma 1 in \cite{2019-Hou-p3} obtained by Hou et al.
  In Theorem \ref{theorem3.5}, if $W_0=1$, then Theorem \ref{theorem3.5} is Lemma 1 in \cite{2019-Hou-p3}.

 We get the upper and lower bounds for {the spectral radius} of $\mathcal{A} (G)$ in Theorem \ref{theorem3.9}, where $G\in \mathbb{Q}_{k,n}$ {with $k,n \geq 3$}.
   Since the proofs for Theorem \ref{theorem3.9} are  similar to those for Theorem \ref{theorem3.8}, we omit it here.

\begin{theorem}\label{theorem3.9} (\textbf{The bounds for {the spectral radius} of $\mathcal{A} (G)$})
  Let  $G\in \mathbb{Q}_{k,n}$, where $k,n \geq 3$.
{Then $ \delta\leq \rho(G) \leq \alpha.$}
\end{theorem}

 \begin{theorem}\label{theorem4.5} (\textbf{$H^+$-eigenvalue of $\mathcal{A} (G)$})
  Let  $G\in \mathbb{Q}_{k,n}$, where $k,n \geq 3$.
 Then zero is an $H^+$-eigenvalue of $\mathcal{A} (G)$ and
 \textcolor{black}{$\bm{e}^{(i)}$ is an $H^+$-eigenvector}
   of $\mathcal{A} (G)$ corresponding to zero, where $i\in[n]$.
\end{theorem}

\noindent\textbf{Proof}.
\textcolor{black}{For any $i,j\in [n]$, we get
\begin{align}
 \left(\mathcal{A} (G) \bm{e}^{(i)}\right)_j
 &= \sum\limits_{j_2 ,\ldots, j_k=1}^{n} a_{j j_2 \ldots j_{k}}e^{(i)}_{j_2} \cdots e^{(i)}_{j_k} \nonumber\\
  &= a_{i i \ldots i}e^{(i)}_i \cdots e^{(i)}_i = 0 = 0 \cdot e^{(i)}_j.  \label{6.3}
     \end{align}}
  It is noted that  (\ref{6.3}) follows from (\ref{A}).
  Thus, zero is an $H^+$-eigenvalue of $\mathcal{A} (G)$ and $\bm{e}^{(1)}, \ldots, \bm{e}^{(n)}$ are the  $H^+$-eigenvectors of $\mathcal{A} (G)$ corresponding to zero.
~~$\Box$

\begin{acknowledgments}
\noindent\textbf{Acknowledgments}\\
 The work was supported by the Natural Science Foundation of Shanghai under the grant number 21ZR1423500 and the National Natural Science Foundation of China under the grant number 11871040. The authors are indebted to the reviewer's helpful comments.
\end{acknowledgments}

%\bibliographystyle{elsevier_citation_order}  %% Numbered references with article and chapter titles, listed in order of citation
%\bibliography{bibliography20211226wwh-sr}

\end{document}